\documentclass{amsart}
\usepackage{amscd}

\input{tcilatex}
\theoremstyle{definition}
\theoremstyle{remark}
\numberwithin{equation}{section}

\input tcilatex

\begin{document}
\title[The theorem of Mather...]{The theorem of Mather on generic projections for singular varieties}
\author{A. Alzati - E. Ballico - G. Ottaviani}
\maketitle

\begin{abstract}
The theorem of Mather on generic projections of smooth algebraic varieties
is also proved for the singular ones.
\keywords{jet spaces, generic projections.}
\subjclass{primary 32C40; secondary 14B05, 14N05.}
\end{abstract}

\smallskip

\section{Introduction}

In [1] it appeared a self-contained proof of the following transversality
theorem of Mather on generic projections (see [2]) in the setting of
Algebraic Geometry:

\medskip

\textbf{Theorem (1,1):} let $X$ be a smooth subvariety of codimension $c$ of
the complex projective space \textbf{P}$^{n}$ . Let $T$ be any linear
subspace of \textbf{P}$^{n}$ of dimension $t$ such that $T\cap X$ $=$ $%
\emptyset $ (so $t\leq c-1$ ). For any $i_{1}\leq t+1$ let $X_{i_{1}}$ $=$ $%
\{x\in X|\dim [(TX)_{x}\cap T]=i_{1}-1\}$ (the dimension of $\emptyset $ is $%
-1$). When $X_{i_{1}}$ is smooth, for any $i_{2}\leq i_{1}$ define $%
X_{i_{1}},_{i_{2}}=$ $\{x\in X_{i_{1}}|\dim [(TX_{i_{1}})_{x}\cap
T]=i_{2}-1\}$ and so on; for $i_{k\text{ }}\leq $ $...$ $i_{2}\leq i_{1}$
define (when possible) $X_{i_{1}},.....,_{i_{k}}$. For $T$ in a Zariski open
set of the Grassmannian $Gr(\mathbf{P}^{t},\mathbf{P}^{n}),$ each $%
X_{i_{1}},.....,_{i_{k}}$ is smooth (and so the above definitions are
possible) until (increasing $k$) it becomes empty and its codimension $%
\upsilon _{I}$ in $X$ can be calculated (where $I=(i_{1},i_{2}.....,i_{k})$).

\medskip

We refer to [1] for the calculation of $\upsilon _{I}$ and for comments and
remarks about the theorem.

This theorem was stated for smooth subvarieties of \textbf{P}$^{n}$ but the
same proof can also be used for the smooth open set $X$ of a singular
algebraic variety $Y$ except for the crucial th. 3.15, (p. 409 of [1]), in
which the compactness of $X$ is needed.

In this short note we want to replace the proof in [1] with a little longer
proof which works also in the case under examination. We obtain the
following:

\medskip

\textbf{Theorem (1,2):} theorem (1,1) still holds if $X$ is replaced with
the smooth open subvariety of a possibly singular projective variety $Y$.

\section{Background}

Let $Y$ be a singular algebraic subvariety of the n-dimensional projective
space \textbf{P}$^{n}$ over the complex numbers. Let $X$ be the smooth open
set of $Y$. First of all we outline the proof of Mather's theorem given in
[1] and we introduce some notation.

Fix an integer $t$ with $0\leq t\leq c-1$. Let $L$ be a $(n-t-1)$%
-dimensional linear subspace of \textbf{P}$^{n}$.

Let $F=$ $\{\mathbf{P}^{t}\in Gr(\mathbf{P}^{t},\mathbf{P}^{n})|\mathbf{P}
^{t}\cap X=\emptyset $ and $\mathbf{P}^{t}\cap L=\emptyset \}$. For any $%
f\in F$ let $p_{f}:X$ $\rightarrow L$ be the linear projection centered in $%
f $ and let $j^{k}p_{f}$ be its $k$-jet ($j^{k}p_{f}:X$ $\rightarrow
J^{k}(X,L) $ sends every $x\in X$ into the $k$-jet of $p_{f}$ in $x$, see
[1] for the definition of $J^{k}(X,L)$ ). Let $I=(i_{1},i_{2}.....,i_{k})$
be any sequence of integers with $(i_{1}\geq i_{2}.....\geq i_{k\text{ }%
}\geq $ $0)$ .

Let $g:X\times $ $F\rightarrow J^{k}(X,L)$ be given by: $g(x,f)=$ $%
(j^{k}p_{f})_{x}$.

The proof of Mather's theorem is divided into two steps:

1) define in $J^{k}(X,L)$ some submanifolds $\Sigma ^{I}$ with the property
that $j^{k}p_{f}^{-1}$ $(\Sigma ^{I})=X_{I\text{ }}$(when $X_{I\text{ }}$are
defined), this definition is not trivial and it is due to Boardman: $\Sigma
^{I}$ are the so called Thom-Boardman singularities, they are smooth,
locally closed and of codimension $\upsilon _{I}$;

2) show that there exists a Zariski open set $U\in F$ such that for any $%
f\in U,$ $j^{k}p_{f}:X$ $\rightarrow J^{k}(X,L)$ is transversal to $\Sigma
^{I}$.

The proof of step 1) runs exactly as in [1].

To prove step 2) firstly we remark, (see [1], prop. 3.13), that for any
smooth subvariety $W\subset J^{k}(X,L)$ there exists a Zariski open set $%
U\in F$ such that for any $f\in U,$ $j^{k}p_{f}:X$ $\rightarrow J^{k}(X,L)$
is transversal to $W$ if $g$ is transversal to $W$. Secondly we give the
following definition. Let $\varphi :X$ $\rightarrow J^{k}(X,L)$ be a
holomorphic map and let $W\subset J^{k}(X,L)$ be a smooth subvariety, then
define:

$\delta (\varphi ,W,x)=0$ if $\varphi (x)\notin W$

$\delta (\varphi ,W,x)=\dim [J^{k}(X,L)]-\dim [TW_{\varphi (x)}+d\varphi
(TX)_{x}]$ if $\varphi (x)\in W$ where $TW$ and $TX$ are the tangent spaces
and $d$ stands for the usual differential. Note that $\delta (\varphi ,W,x)%
\geq 0$ and that $\varphi $ is transversal to $W$ at $x$ if and only if $%
\delta (\varphi ,W,x)=0$ .

As in [1], th. 3.10 and 3.11, it can be shown that for $W=\Sigma ^{I}\subset
J^{k}(X,L)$ the following condition $(*)$ is satisfied:

$(*)$ $\delta (g,W,(x,f))\leq \delta (j^{k}p_{f},W,x)$ for any $(x,f)\in X%
\times $ $F$ and 

equality holds if and only if $\delta (j^{k}p_{f},W,x)=0$.

Therefore to prove step 2) all that we need is the following:

\medskip

\textbf{Theorem (2,1): }with the previous notation, assume that condition $%
(*)$ is satisfied for some smooth subvariety $W\subset J^{k}(X,L)$; then
there exists a Zariski open set $U\in F$ such that for any $f\in U,$ $%
j^{k}p_{f}:X$ $\rightarrow J^{k}(X,L)$ is transversal to $W$.

\medskip

The proof of this theorem (th. 3.15 in [1]) must be rewritten in our case.
In \S 3 we will give this proof and so we will also prove theorem (1.2).

\section{Proof of theorem (2.1)}

Let us define $\delta _{g}=Sup_{(x,f)\in X\times F}\{\delta (g,W,(x,f))\}$,
moreover let us define $A=\{(x,f)\in X\times $ $F|\delta (g,W,(x,f))=\delta
_{g}\}\subset X\times F$, $A$ is a Zariski closed set in $X\times $ $F$.
Note that if $\delta _{g}=0$ theorem (2.1) is true (see th.\ 3.13 in [1]),
so we can assume $\delta _{g}\neq 0$ and $A\neq \emptyset $.

Let $\pi _{2}:X\times F\rightarrow F$ be the natural projection. $X\times F$
is equipped with the induced Zariski topology from $Y\times $ $F$. Let $%
\overline{A}$ be the Zariski closure of $A$ in $Y\times $ $F$; let $\pi %
_{3}:Y\times $ $F\rightarrow F$ be the natural projection, $\pi _{3}(%
\overline{A})$ is a Zariski closed set of $F$. If $\pi _{3}(\overline{A})$
is a proper subset of $F,$ we can consider $F^{\prime }$ $=F\backslash \pi %
_{3}(\overline{A})$ and $g^{\prime }$ $=g_{|X\times F^{\prime }}$. The
assumptions of the theorem are true for $F^{\prime }$ and $g^{\prime }$ and $%
\delta _{g^{\prime }}<$ $\delta _{g}$. If the corresponding $\pi _{3}(%
\overline{A})$ were a proper subset of $F^{\prime }$ we would get $F^{\prime
\prime }$ and $g^{\prime \prime }$ and so on. After a finite number of steps
we would get $F^{\prime }$ and $g^{\prime }$, for which the assumptions
would be still true, with $\delta _{g^{^{\prime }}}$ $=0$, so the theorem
would be proved.

Hence we have only to prove that $\pi _{3}(\overline{A})$ is a proper subset
of $F$.

By contradiction let us assume that $\pi _{3}(\overline{A})$ $=$ $F$, then $%
F=$ $\overline{\pi _{2}(A)}$.

We can choose $(x_{0},f_{0})\in A$ and $z_{0}=(j^{k}g)_{(x_{0},f_{0})}\in W$
. As $\delta (g,W,(x_{0},f_{0}))$ is strictly positive, by assumption we get
that $\delta (j^{k}p_{f_{0}},W,x_{0})$ is strictly positive too, hence $%
j^{k}p_{f_{0}}$ is not transversal to $W$ at $x_{0}$.

$W$ is smooth at $x_{0}$ so it is a local complete intersection, then it is
possible (see [1], proof of th. 3.15) to get a smooth subvariety $W^{\prime
}\subset J^{k}(X,L)$ and a smooth dense open Zariski set $Z\subset X\times F$
such that: $W\subseteq W^{\prime }$, $\dim (W^{\prime })-\dim (W)=\delta _{g}
$, $g$ is transversal to $W^{\prime }$ at $(x,f)$ for any $(x,f)\in Z$.

The holomorphic map $g_{|Z}:Z$ $\rightarrow J^{k}(X,L)$ is transversal to $%
W^{\prime }$ so that $g_{|Z}^{-1}(W^{\prime })=$ $g^{-1}(W^{\prime })\cap Z$
is smooth in $X\times F$.

Let us consider $\pi =\pi _{2_{|g^{-1}(W^{\prime })\cap Z}}=\pi
_{3_{|g^{-1}(W^{\prime })\cap Z}}:g^{-1}(W^{\prime })\cap Z\rightarrow F$.

It is easy to see that:

$(1)$ $\overline{\pi _{2}(A\cap Z)}=F$

hence $F=$ $\overline{\pi _{2}(Z)}$. Moreover $F=$ $\overline{\pi %
_{2}(g^{-1}(W^{\prime }))}$ otherwise there would exist a Zariski open set $%
B\subset F$ such that $B\cap $ $\overline{\pi _{2}(g^{-1}(W^{\prime }))}%
=\emptyset $, hence for any $f\in B$ and for any $x\in X$, $(x,f)\notin
g^{-1}(W^{\prime })$, i.e. $g(x,f)\notin W^{\prime }$, i.e. $g(x,f)\notin W$%
, i.e. for any $f\in B$ and for any $x\in X$, $\delta (g,W,(x,f))=0$ and the
theorem would be immediately proved (see th. 3.13 of [1]).

It follows: $\overline{\pi _{2}(g^{-1}(W^{\prime })\cap Z)}$ $\subseteq $ $%
\overline{\pi _{2}(g^{-1}(W^{\prime }))\cap \pi _{2}(Z)}$ $\subseteq $ $%
\overline{\pi _{2}(g^{-1}(W^{\prime })})$ $\cap \overline{\pi _{2}(Z)}=$ $F$%
, therefore:

$(2)$ $\overline{\pi (g^{-1}(W^{\prime })\cap Z)}=F$.

Now we consider the holomorphic map $\pi :g^{-1}(W^{\prime })\cap Z$ $%
\rightarrow F$ between smooth manifolds, as $(2)$ holds there exists a
Zariski open set $D\subset F$ such that for any $f\in D$ $\pi ^{-1}(f)$ is
smooth and of the expected codimension.

By $(1)$ $[\pi _{2}(A\cap Z)]\cap D\neq \emptyset $, then we can choose $%
f_{1}\in [\pi _{2}(A\cap Z)]\cap D$ such that $\pi ^{-1}(f_{1})$ is smooth,
of the expected codimension and biholomorphic to a Zariski open set of $%
(j^{k}p_{f_{1}})^{-1}(W^{\prime })\subset X$. We can also choose $x_{1}\in X$
such that $(j^{k}p_{f_{1}})^{-1}(W^{\prime })$ is smooth, of the expected
codimension and smooth at $x_{1}$. This fact implies that $j^{k}p_{f_{1}}$
is transversal to $W^{\prime }$ at $x_{1}$, (see [1], th. 1.2), i.e. $\delta
(j^{k}p_{f_{1}},W^{\prime },x_{1})=0$.

On the other hand $f_{1}\in \pi _{2}(A\cap Z)$, hence it is possible to
choose $x_{1}\in X$ such that $(x_{1},f_{1})\in A$, i.e. $\delta
(g,W,(x_{1},f_{1}))=$ $\delta _{g}$.

Let $z_{1}=(j^{k}p_{f_{1}})_{x_{1}}$ then:

$\delta (j^{k}p_{f_{1}},W^{\prime },x_{1})=\dim [J^{k}(X,L)]-\dim
[(TW^{\prime })_{z_{1}}+dj^{k}p_{f_{1}}(TX)_{x_{1}}]$

$\delta (j^{k}p_{f_{1}},W,x_{1})=\dim [J^{k}(X,L)]-\dim
[(TW)_{z_{1}}+dj^{k}p_{f_{1}}(TX)_{x_{1}}]$

and $0=\delta (j^{k}p_{f_{1}},W^{\prime },x_{1})\geq \delta
(j^{k}p_{f_{1}},W,x_{1})-\delta _{g}$.

But assumption $(*)$ and the fact that $(x_{1},f_{1})\in A$ imply:

$0\geq \delta (j^{k}p_{f_{1}},W,x_{1})-\delta _{g}>\delta
(g,W,(x_{1},f_{1}))-\delta _{g}=\delta _{g}-\delta _{g}=0$, contradiction!

\section{Cones}

In this brief section we want to remark that when $Y$ is a cone it is
possible to use Mather's theorem (1.1).

For instance let us assume that $Y$ is a cone in \textbf{P}$^{n}$ of vertex $%
V$ on a smooth subvariety $B$ of \textbf{P}$^{n}$ whose span is \textbf{P}$%
^{s}$ with $\dim (Y)=\gamma =b+v+1,\dim (B)=b,\dim (V)=v,n=s+v+1$.

Let $T$ be a generic $t$-dimensional subspace of \textbf{P}$^{n}$ with: $%
T\cap Y=\emptyset $, $t\leq \gamma -1$, $\gamma \geq (t+1)(n-\gamma )$. Let $%
Y_{t+1}=\{y\in Y|$ $y$ is a smooth point, $(TY)_{y}\supset T\}$.

If $Y$ were smooth Mather's theorem (1.1) would say that, for generic $T$, $%
Y_{t+1}$ is a smooth subvariety of $Y$ and $\dim (Y_{t+1})=\gamma
-(t+1)(n-\gamma )$, in our case we have:

\medskip

\textbf{Proposition:} the closure of $Y_{t+1}$ is a cone of dimension $%
\gamma -(t+1)(n-\gamma )$ with vertex $V$ over a smooth variety.

\medskip

As $Y$ is a cone we remark that $t\leq s-1$ ($t\leq n-\gamma -1=s-b-1$),
hence there exists a linear subspace $H\simeq $ \textbf{P}$^{s}$ in \textbf{P%
}$^{n}$ such that $H\supset T$ and $H\cap V=\emptyset $. We can assume that $%
B=H\cap Y$ and we can apply th. (1.1) to \textbf{P}$^{s}$, $T$ and $B$ as $B$
is smooth, $T\cap B=\emptyset $ and $T$ is generic in \textbf{P}$^{s}$ with
respect to $B$. If $t\leq b-1$ and $b\geq (t+1)(s-b)$ (for instance when $%
t=0 $ and $2b\geq s$) then $B_{t+1}=\{y\in B|(TB)_{y}\supset T\}$ is a
smooth subvariety of $B$ and $\dim (B_{t+1})=b-(t+1)(s-b)$. On the other
hand $(TB)_{y}\supset T$ if and only if $(TY)_{y}\supset T$ as $%
(TY)_{y}=\langle V,(TB)_{y}\rangle $ i.e. $(TB)_{y}=(TY)_{y}\cap H$, hence $%
Y_{t+1}\cap H=B_{t+1}$ and the closure in $Y$ of $Y_{t+1}$ is another cone
of vertex $V$ over $B_{t+1}$. This cone has dimension $b-(t+1)(s-b)+v+1=%
\gamma -(t+1)(n-\gamma )$ which is exactly the expected dimension when $Y$
is smooth.

\section{References}

[1] A. Alzati, G. Ottaviani: ''The theorem of Mather on generic projections
in the setting of Algebraic Geometry'' Manuscr. Math. \textbf{74} 391-412
(1992).

[2] J.\ N.\ Mather: ''Generic projections'' Ann. of Math. \textbf{98}
226-245 (1973).

\bigskip

\bigskip

\bigskip

{\footnotesize All authors are members of Italian GNSAGA. Work supported by
Murst funds.}

{\footnotesize Addresses:}

{\footnotesize A.\ Alzati: Dip. di Matematica Univ. di Milano, via C.\
Saldini 50 20133-Milano (Italy). }

{\footnotesize E-mail: alzati@mat.unimi.it}

{\footnotesize E.\ Ballico: Dip. di Matematica Univ. di Trento, via
Sommarive 14 38050-Povo Trento (Italy). }

{\footnotesize E-mail: ballico@science.unitn.it}

{\footnotesize G.\ Ottaviani: Dip di Matematica Univ. di Firenze, viale
Morgagni 67/A 50134-Firenze (Italy). }

{\footnotesize E-mail: ottavian@udini.math.unifi.it}

\end{document}